\documentclass{amsart}
\usepackage{amssymb,amsmath,epsfig}
\usepackage[british]{babel}

\newtheorem{theorem}{Theorem}[section]
\newtheorem{proposition}[theorem]{Proposition}
\newtheorem{lemma}[theorem]{Lemma}
\newtheorem{corollary}[theorem]{Corollary}
\newtheorem{letterthm}{Theorem}

\theoremstyle{definition}
\newtheorem{definition}[theorem]{Definition}

\newcommand{\R}{\mathbb{R}}
\newcommand{\Z}{\mathbb{Z}}
\newcommand{\me}{\stackrel{\mathrm{ME}}{\sim}}

\title[Limit groups and measure equivalence]{Limit groups, positive-genus
towers and measure equivalence}
\date{6th November 2006}

\author{Martin R.~Bridson}
\address{Dept.\ of Mathematics, Imperial College, London SW7 2AZ, UK.}
\email{m.bridson@imperial.ac.uk}
\thanks{}

\author{Michael Tweedale}
\address{Dept.\ of Mathematics, University of Bristol, Bristol BS8 1TW, UK.}
\email{m.tweedale@bristol.ac.uk}
\thanks{}

\author{Henry Wilton}
\address{Dept.\ of Mathematics, 1 University Station C1200, Austin, TX 78712,
USA.}
\email{henry.wilton@math.utexas.edu}
\thanks{This work was supported in part by grants from the EPSRC. The first
author is also supported by a Royal Society Wolfson Research Merit
Award.}

\begin{document}
\begin{abstract}
By definition, an $\omega$-residually free tower is positive-genus if all
surfaces used in its construction are of positive genus.  We prove that every
limit group is virtually a subgroup of a positive-genus $\omega$-residually
free tower. By combining this construction with results of Gaboriau, we prove
that elementarily free groups are measure equivalent to free groups.
\end{abstract}

\maketitle\let\languagename\relax\sloppy

\section{Introduction}
\noindent Measure equivalence was introduced by M.~Gromov in~\cite{Gro93} as a
measure-theoretic analogue of quasi-isometry.  After proving that two finitely
generated groups are quasi-isometric if and only if they admit commuting,
properly discontinuous actions by homeomorphisms on some locally compact space
with \emph{compact} fundamental domain, he went on to define measure
equivalence by requiring instead that the groups admit commuting, free,
measure-preserving actions on a Lebesgue measure space with \emph{finite
measure} fundamental domain. The \emph{uniform} lattices in a fixed locally
compact second countable group are all quasi-isometric, whereas \emph{all}
lattices in such a group are measure equivalent, regardless of whether they
are uniform or not.  For motivation and background, we refer the reader to the
papers of Damien Gaboriau, particularly \cite{Gab00}.

Much progress has been made in distinguishing measure-equivalence classes
(see, for example, \cite{Adams95}, \cite{Furman99a}, \cite{Furman99b},
\cite{Gab00}, \cite{Gab02} and~\cite{Zim84}), but there have been many fewer
constructions of examples of measure equivalent groups.  The only groups whose
measure-equivalence classes are completely classified are finite groups,
amenable groups, and lattices in simple Lie groups of higher rank
(see~\cite{Furman99a}). In particular, the measure-equivalence class of
non-abelian free groups is still quite poorly understood: there are some
obvious examples like the fundamental groups of compact surfaces of genus
$\ge2$, which are also lattices in $\mathrm{SL}(2,\R)$, but beyond this we
know very little.

In~\cite{Gab05}, Gaboriau constructs some new examples of groups measure
equivalent to free groups, encapsulated in the following theorem.

\begin{theorem}[Gaboriau]\label{Amalgam of handles}
Let $\Sigma$ be a compact orientable surface of positive genus, with one
boundary component. Then the iterated amalgamated product
\[
\pi_1(\Sigma)*_{\langle\partial\Sigma\rangle}\pi_1(\Sigma)*_{\langle\partial\Sigma\rangle}\ldots*_{\langle\partial\Sigma\rangle}\pi_1(\Sigma)
\]
is measure equivalent to a free group.
\end{theorem}
The amalgamated product in the above statement is an example of a \emph{limit
group} (indeed, it is an elementarily free group).  Gaboriau \cite{Gab05}
asked if all limit groups are measure equivalent to free groups.

Limit groups have been studied under a variety of names and guises:
see~\cite{Baum67}, \cite{CG04}, \cite{KM98a} and~\cite{KM98b}. The name limit
group was introduced by Z.~Sela in his solution of the Tarski problem
(see~\cite{Se1}~\emph{et~seq.}; cf.~\cite{KM98a},
\cite{KM98b}~\emph{et~seq.}).  The \emph{elementary theory} of a group $G$ is
the set of first-order sentences that are true in $G$.  The Tarski problem
asks which groups have the same elementary theory as the free group of
rank~$2$.  The \emph{existential theory} consists of those sentences involving
only existential quantifiers.  Limit groups turn out to be precisely the
finitely generated groups with the same existential theory as a free group.
Another class, still more closely related to free groups, is the class of
\emph{elementarily free groups}: those groups with the same elementary theory
as a free group.  It turns out that every elementarily free group is a
(Gromov) hyperbolic limit group.  However, there are hyperbolic limit groups
that are \emph{not} elementarily free.

Using results of Kharlampovich--Myasnikov and Sela~\cite{KM98a, KM98b, Se1},
one can give more constructive definitions of limit groups and elementarily
free groups.  For the purposes of this paper, both classes are defined in
terms of $\omega$-residually free towers, which are in turn defined as the
fundamental groups of certain complexes, inductively constructed from graphs,
surfaces and tori. The limit groups are precisely the finitely generated
subgroups of $\omega$-residually free towers, but not every limit group is
itself such a tower (see section~\ref{Limit groups}).  We call a tower
\emph{positive-genus} if every surface used in its construction is of positive
genus. Our main result concerning towers is the following.

\begin{letterthm}[Theorem~\ref{EFGs are virtually
  positive-genus}]\label{Technical result}
  Every limit group $L$ is virtually a subgroup of a positive-genus
  $\omega$-residually free tower. If $L$ is elementarily free, then the tower
  can be chosen to be elementarily free (hyperbolic).
\end{letterthm}

In general, it is necessary to pass to a finite-index subgroup in order to
embed a limit group in a positive-genus tower.  For example, if $L$ is
constructed by amalgamating the fundamental group of a punctured sphere along
its boundary components with another limit group $L'$, in such a way that the
resulting two-vertex graph-of-groups decomposition of $L$ is the cyclic JSJ
decomposition, then $L$ does not embed in a positive-genus tower.

Using theorem \ref{Technical result} and results of~\cite{Gab05}, we deduce
the following theorem, which partially answers Gaboriau's question.

\begin{letterthm}[Theorem~\ref{EFGs are ME to free}]
  Every elementarily free group is measure equivalent to a free group.
\end{letterthm}

This paper is organized as follows.  In section~1 we construct various
finite-index subgroups of towers and prove theorem~\ref{Technical result}. In
section~2 we recapitulate some useful results of~\cite{Gab05} and prove that
elementarily free groups are measure equivalent to free groups. In section~3
we discuss methods of attacking the case of all limit groups.

We thank the referee for his or her helpful comments.

\section{Limit groups}\label{Limit groups}
\subsection{$\omega$-residually free towers}
For an introduction to the theory of limit groups see, for instance,
\cite{BF03} or~\cite{CG04}.

\begin{definition}\label{Tower}
  An \emph{$\omega$-residually free tower space of height $0$}, denoted $X_0$,
  is a one-point union of finitely many compact graphs, $n$-tori, and closed
  hyperbolic surfaces of Euler characteristic less than~$-1$.

  An \emph{$\omega$-residually free tower space of height $h$}, denoted $X_h$,
  is obtained from an $\omega$-residually free tower space $X_{h-1}$ of height
  $h-1$ by attaching one of two sorts of blocks.
  \begin{enumerate}
    \item \textbf{Quadratic block.}  Let $\Sigma$ be a connected compact
      hyperbolic surface with boundary, with each component either a punctured
      torus or having $\chi\leq-2$.  Then $X_h$ is the quotient of
      $X_{h-1}\sqcup\Sigma$ obtained by identifying the boundary components of
      $\Sigma$ with homotopically non-trivial curves on $X_{h-1}$, in such a
      way that there exists a retraction $\rho:X_h\rightarrow X_{h-1}$. The
      retraction is also required to satisfy the property that
      $\rho_*(\pi_1(\Sigma))$ be non-abelian.

    \item \textbf{Abelian block.} Let $T$ be an $n$-torus, and fix a
      coordinate circle $\gamma$.  Fix a loop $c$ in $X_{h-1}$ that generates
      a maximal abelian subgroup in $\pi_1(X_{h-1})$.  Then $X_h$ is the
      quotient of $X_h\sqcup (S^1\times[0,1])\sqcup T$ obtained by identifying
      $S^1\times\{0\}$ with $c$, and $S^1\times\{1\}$ with $\gamma$.
  \end{enumerate}
  Attaching an abelian block is also often called \emph{extending a
  centralizer}.

  An $\omega$-residually free tower space is called \emph{hyperbolic} if no
  tori are used in its construction; that is, there are no tori in $X_0$ and
  no abelian blocks.
\end{definition}

\begin{definition}
  An \emph{($\omega$-residually free) tower of height $h$}, denoted $L_h$, is
  the fundamental group of an $\omega$-residually free tower space of height
  $h$.
\end{definition}

The following deep theorem from~\cite{Se6} (cf.~\cite{KM98a} \emph{et seq.})
will, for our purposes, serve as a definition of elementarily free groups.

\begin{theorem}
  A group is \emph{elementarily free} if and only if it is the fundamental
  group of a hyperbolic $\omega$-residually free tower space.
\end{theorem}

Towers are examples of limit groups. Another theorem of Sela~\cite{Se2} and
O.~Kharlampovich and A.~Myasnikov~\cite{KM98b} will serve as a definition of
limit groups.

\begin{theorem}\label{Limit groups embed in towers}
  A group is a \emph{limit group} if and only if it is a finitely generated
  subgroup of an $\omega$-residually free tower.
\end{theorem}

A little care is needed here.  While every limit group embeds in a tower,
there exist hyperbolic limit groups that do not embed in hyperbolic towers.

Improving on theorem~\ref{Limit groups embed in towers}, Kharlampovich and
Myasnikov proved that every limit group embeds in a tower in which no surface
pieces are used at level~$0$ and no quadratic blocks are used (\cite{KM98b};
cf.~\cite{CG04}).

We will need the following weak form of the well known fact that limit groups
are precisely the finitely generated $\omega$-residually free groups.

\begin{lemma}\label{Residually free}
  Limit groups are \emph{residually free}; that is, if $L$ is a limit group
  and $g\in L-\{1\}$ then there exists a homomorphism to a free group
  $f:L\rightarrow F$ with $f(g)\neq 1$.
\end{lemma}

A key feature of the definition of a tower is the retraction
$\rho:X_h\rightarrow X_{h-1}$.  If the block at height $h$ is abelian, the
retraction simply projects $T$ onto the coordinate circle $\gamma$, and thence
to $c$.  In both cases, $\rho$ induces a retraction $\rho_*:L_h\rightarrow
L_{h-1}$ on the level of fundamental groups.

An $\omega$-residually free tower space $X_h$ has a natural
graph-of-spaces\footnote{By convention, our graphs of spaces are connected and
have connected vertex and edge spaces.} decomposition $\Gamma_X$, with two
vertex spaces, namely $X_{h-1}$ and the block at height $h$; the edge spaces
are circles.  We will often use the retraction to pull finite covers back from
$X_{h-1}$ to $X_h$. It is worth noting that such pullbacks inherit a similar
graph-of-spaces decomposition from $X_h$.

\begin{lemma}\label{Pullback covers}
  Let $X$ be a CW-complex with a graph-of-spaces decomposition $\Gamma_X$,
  such that there is a retraction $\rho:X\rightarrow X'$ to a vertex space.
  Let $Y'\rightarrow X'$ be a connected covering of degree $d$, and let
  $Y\rightarrow X$ be the connected covering obtained by pulling back along
  $\rho$; that is, $\pi_1(Y)=\rho_*^{-1}(\pi_1(Y'))$.  Then:
  \begin{enumerate}
    \item $Y\rightarrow X$ is of degree $d$ and inherits a graph-of-spaces
      decomposition $\Gamma_Y$; \item the pre-image of $X'$ in $Y$ is a
      (connected) vertex space of $\Gamma_Y$ naturally homeomorphic to $Y'$;
    \item $Y\rightarrow X$ extends $Y'\rightarrow X'$, and $Y$ inherits a
      retraction to $Y'$ covering $\rho$.
  \end{enumerate}
\end{lemma}

A tower $L_h$ inherits, by the Seifert--van Kampen Theorem, a graph-of-groups
decomposition $\Gamma_L$  from the graph-of-spaces decomposition $\Gamma_X$ of
the associated $\omega$-residually free tower space $X_h$.

\subsection{Positive-genus towers}
The purpose of this section is to prove that, up to finite index, the
quadratic blocks in definition \ref{Tower} can be assumed to have positive
genus.  A compact, connected surface $\Sigma$ with Euler characteristic
$\chi(\Sigma)$ and $b(\Sigma)$ boundary components is of \emph{positive genus}
if $\chi(\Sigma)+b(\Sigma)\leq 0$. Note that, in particular, all finite covers
of such $\Sigma$ also have positive genus.

\begin{definition}
  An $\omega$-residually free tower space is \emph{positive-genus} if every
  quadratic block used in its construction is of positive genus.  A tower is
  \emph{positive-genus} if it is the fundamental group of a positive genus
  $\omega$-residually free tower space.
\end{definition}

Our strategy for obtaining positive-genus quadratic blocks is to identify
connected $p$-sheeted coverings (here $p$ is a prime number) that restrict to
a $p$-sheeted covering on each boundary component.  We achieve this by passing
to a finite-index subgroup of the tower that admits a map to $\Z/p\Z$ which
maps each attaching loop of the top quadratic block non-trivially.  In
particular, we must arrange for the attaching loops to become non-trivial in
homology.

\begin{lemma}\label{Curve lifting}
  If $G$ is a finitely generated residually free group and $g_1,\ldots,g_m$ is
  a finite collection of elements of $G$, then there exists a finite-index
  subgroup $K\subset G$ so that, whenever $k\in K$ is a power of a conjugate
  of some $g_i$, $k$ has infinite order in $H_1(K)$.
\end{lemma}
\begin{proof}
  Let $F$ denote the free group of rank~$2$. Since $L$ is residually free, for
  each $i$ there exists a homomorphism $f_i:L\rightarrow F$ with $f_i(g_i)\neq
  1$.  By M.~Hall's theorem~\cite{Hall49}, there exists a finite-index
  subgroup $F_i\subset F$ containing $f_i(g_i)$, such that $f_i(g_i)$ is
  primitive in $H_1(F_i)$.  Let $K$ be the finite-index subgroup $\bigcap_i
  f_i^{-1}(F_i)$ of $G$.  If $k\in K$ is a power of a conjugate of $g_i$ then,
  since $f_i(g_i)$ has infinite order in $H_1(F_i)$, it follows that $k$ has
  infinite order in $H_1(K)$.
\end{proof}

We can now construct a map to $\Z/p\Z$ as required.

\begin{lemma}\label{The map to Z/p}
  Let $K$ be a finitely generated group and let $k_1,\ldots,k_n$ be a
  collection of elements in $K$ that are all of infinite order in $H_1(K)$.
  Then, for all sufficiently large primes $p$, there exists a homomorphism
  $\varphi:G\rightarrow \Z/p\Z$ so that $\varphi(k_j)$ is non-trivial for all
  $j$.
\end{lemma}
\begin{proof}
  Since $\Z^n$ is $\omega$-residually free, there exists a homomorphism
  $K\rightarrow\Z$ under which each $k_j$ has non-trivial image.  Choose a
  prime $p$ that does not divide any of the images of the $k_j$ in $\Z$.  In
  particular, each $k_j$ has non-trivial image under the composition
  \[
  \varphi:K\rightarrow H_1(K)\rightarrow\Z\rightarrow\Z/p\Z.\qedhere
  \]
\end{proof}

We shall apply the preceding lemmas to the height $h-1$ subspace $X_{h-1}$ and
pull back to the full tower to obtain the positive-genus cover that we seek.

Recall that, for $X$ a topological space, $c:S^1\rightarrow X$ a loop, and
$Y\rightarrow X$ a covering map, the \emph{elevations} of $c$ to $Y$ are the
minimal connected covers $\hat{S}^1\rightarrow S^1$ such that
$\hat{S}^1\rightarrow X$ lifts to $Y$.  Fixing basepoints, it follows from
standard covering space theory that $\pi_1(\hat{S}^1)$ is the pre-image of
$\pi_1(Y)$ in $\pi_1(S^1)$.

\begin{proposition}\label{Positive-genus cover}
  Let $X_h$ be an $\omega$-residually free tower space, constructed by
  attaching a quadratic block $\Sigma$ to a space $X_{h-1}$ of height $h-1$.
  Then there exists a connected cover $Z_h\rightarrow X_h$ with an inherited
  graph-of-spaces decomposition $\Gamma_Z$, with one vertex space a connected
  cover $Z_{h-1}\rightarrow X_{h-1}$, and the remaining vertex spaces
  connected covers $\bar{\Sigma}_i\rightarrow \Sigma$, where each
  $\bar{\Sigma}_i$ has positive genus.  The retraction $\rho:X_h\rightarrow
  X_{h-1}$ pulls back to a retraction $Z_h\rightarrow Z_{h-1}$.
\end{proposition}
\begin{proof}
  Let $c_1,\ldots,c_m$ be the images of the boundary curves of $\Sigma$ in
  $X_{h-1}$.  Since $\pi_1(X_{h-1})$ is residually free, the finite-index
  subgroup provided by lemma~\ref{Curve lifting} corresponds to a
  finite-sheeted covering $Y_{h-1}\rightarrow X_{h-1}$ so that if
  $d_1,\ldots,d_n$ are the elevations of the $c_i$, each $d_j$ is of infinite
  order in homology.  Let $Y_h\rightarrow X_h$ be the covering obtained by
  pulling back along the retraction $\rho$, with graph-of-spaces decomposition
  $\Gamma_Y$.

  By lemma~\ref{The map to Z/p}, there exists a homomorphism
  $\varphi:\pi_1(Y_{h-1})\rightarrow \Z/p\Z$ with $\varphi(d_j)\neq 0$ for
  each $j$.  Let $Z_{h-1}\rightarrow Y_{h-1}$ be the covering corresponding to
  the kernel of $\varphi$.  Finally, pull the covering $Z_{h-1}\rightarrow
  Y_{h-1}$ back along the retraction $Y_h\rightarrow Y_{h-1}$ to give a
  covering $Z_h\rightarrow Y_h$, with graph-of-spaces decomposition
  $\Gamma_Z$.

  The key point to observe is that each edge space of $\Gamma_Y$ is only
  covered by one edge space of $\Gamma_Z$. Indeed, $Z_h\rightarrow Y_h$ is a
  covering of degree $p$, but the image in $\Z/p\Z$ of every cycle $d_j$ has
  order $p$, so an elevation of it to $Z_h$ covers $d_j$ with degree $p$.
  Thus $d_j$ only has one elevation to $Z_h$.

  It follows that the underlying graphs of $\Gamma_Y$ and $\Gamma_Z$ are the
  same.  Consider a surface vertex $\bar{\Sigma}_i$ of $\Gamma_Z$, covering a
  surface vertex $\Sigma_i$ of $\Gamma_Y$. By construction
  $b(\bar{\Sigma}_i)=b(\Sigma_i)$, so we have
  \[
  \chi(\bar{\Sigma}_i)+b(\bar{\Sigma}_i)=p\ \chi(\Sigma_i)+b(\Sigma_i).
  \]
  Since $\chi(\Sigma_i)\leq -1$ and $\chi(\Sigma_i)+b(\Sigma_i)\leq 2$, it
  follows that $\bar{\Sigma}_i$ has positive genus for $p\geq 3$.
\end{proof}

The above result is all we need to prove that elementarily free groups are
measure equivalent to free groups.  It is perhaps more cleanly expressed,
however, in terms of the following theorem, which we believe to be of
independent interest.

\begin{theorem}\label{EFGs are virtually positive-genus}
  Every limit group $L$ has a finite-index subgroup $M$ that is a subgroup of
  a positive-genus tower $P$.  If $L$ is elementarily free then $P$ can be
  taken to be elementarily free (hyperbolic).
\end{theorem}

By theorem~\ref{Limit groups embed in towers} it suffices to prove the theorem
for towers.  More precisely, we prove the following.

\begin{proposition}
  Let $L_h$ be a tower of height $h$.  Then there exists a finite-index
  subgroup $M_h\subset L_h$ that embeds into a positive-genus tower $P_h$.  If
  $L_h$ is elementarily free then $P_h$ can be taken to be elementarily free.
  If $A\subset M_h$ is a maximal abelian subgroup then $A$ is also maximal
  abelian in $P_h$.
\end{proposition}
\begin{proof}
  The proof is by induction on height.  By definition, every level 0 tower is
  positive-genus.  Consider $L_h$ the fundamental group of an
  $\omega$-residually free tower space $X_h$ of height $h$, obtained as usual
  by attaching a block to a height $h-1$ space $X_{h-1}$ with fundamental
  group $L_{h-1}$.

  First we consider the case of a quadratic block $\Sigma$.  By induction,
  $L_{h-1}$ has a finite-index subgroup $M_{h-1}$ that is a subgroup of a
  positive-genus tower $P_{h-1}$.  By proposition~\ref{Positive-genus cover},
  $L_h$ has a finite-index subgroup $K_h=\pi_1(Z_h)$ with graph-of-groups
  decomposition $\Gamma_K$, with one vertex labelled by $K_{h-1}$ a
  finite-index subgroup of $L_{h-1}$ and the remaining vertices labelled by
  the fundamental groups of surfaces of positive genus, amalgamated with
  $K_{h-1}$ along boundary components.  Set $M_h=\rho_*^{-1}(K_{h-1}\cap
  M_{h-1})$.  Then $M_h$ inherits a graph-of-groups decomposition $\Gamma_M$,
  with one vertex labelled by $K_{h-1}\cap M_{h-1}$ and the remainder by
  fundamental groups of surfaces of positive genus, amalgamated with
  $K_{h-1}\cap M_{h-1}$ along boundary components.  The retraction
  $\rho_*:L_h\rightarrow L_{h-1}$ restricts to a retraction $M_h\rightarrow
  K_{h-1}\cap M_{h-1}$. Enlarge $\Gamma_M$ to $\Gamma_P$ by replacing
  $K_{h-1}\cap M_{h-1}$ with $P_{h-1}$.  Extending $\rho_*$ by the identity on
  $P_{h-1}$, it is clear that $P_h=\pi_1(\Gamma_P)$ is a positive-genus tower.

  The case of an abelian block $T$ is similar.  By induction, there exists a
  finite-index subgroup $M_{h-1}\subset L_{h-1}$ that embeds in a positive
  genus tower $P_{h-1}$.   The pullback $M_h=\rho_*^{-1}(M_{h-1})$ inherits a
  graph-of-groups decomposition $\Gamma_M$, with one vertex labelled by
  $M_{h-1}$ and the remainder by finitely generated free abelian groups.  Each
  abelian vertex has a coordinate factor amalgamated with a cyclic maximal
  abelian subgroup of $M_{h-1}$.  Enlarge $\Gamma_M$ to $\Gamma_P$ by
  replacing $M_{h-1}$ by $P_{h-1}$.  Since cyclic maximal abelian subgroup of
  $M_{h-1}$ are maximal abelian in $P_{h-1}$, the resulting fundamental group
  $P_h=\pi_1(\Gamma_P)$ is a positive-genus tower.

  It remains to show that a maximal abelian subgroup $A$ of $M_h$ is maximal
  abelian in $P_h$.  For this we need a little Bass--Serre Theory.  The result
  follows by induction on height if $A$ is elliptic in $\Gamma_P$, so assume
  that $A$ is not elliptic and let $g\in P_h$ be such that $[g,a]=1$ for all
  $a\in A$.

  Consider the Bass--Serre tree $T_P$ of $\Gamma_P$.  Let $v$ be a vertex of
  $T_P$ stabilized by $P_{h-1}$.  Certain edges $\{e_i\}$ adjoining $v$
  correspond to the cosets of $(M_h:C)\subset (P_h:C)$ for some edge group $C$
  of $\Gamma_P$.  Consider the subspace $T_M\subset T_P$ given by the
  $M$-translates of $v$ and the $e_i$. Then $T_M$ is connected since the
  normal form in $\Gamma_M$ of an $m\in M_h$ specifies a path from $v$ to $mv$
  and furthermore, by construction, $T_M$ is in fact the Bass--Serre tree of
  $\Gamma_M$.

  Let $l\subset T_P$ be the unique line stabilized by $A$.  Fix an edge $e$ in
  $l$. Then $ge$ is an edge of $l$, so lies in $T_M$.  There is only one
  $M_h$-orbit of $e$ in $T_M$, so there exists $m\in M_h$ such that $me=ge$.
  The stabilizer of $e$ lies in $M_h$, so it follows that $g\in M_h$.  Since
  $A$ was maximal abelian in $M_h$, we get $g\in A$.
\end{proof}

\section{Measure equivalence}
\noindent We are now in a position to use the results of~\cite{Gab05} to prove
that elementarily free groups are measure equivalent to free groups.

\subsection{Definition and properties}
\begin{definition}
  Two countable groups $G_1,G_2$ are \emph{measure equivalent} if there exist
  commuting, measure-preserving, (essentially) free actions on some measure
  space $(\Omega,m)$, such that the action of $G_i$ admits a finite measure
  fundamental domain.  Write
  \[
  G_1\me G_2.
  \]
\end{definition}

The standard examples of measure-equivalent groups are commensurable groups
and lattices in the same locally compact second countable group.  We will not
use the definition of measure equivalence directly, but deduce our result from
the following properties.

\begin{theorem}[$\mathbf{P_{ME}7}$ in~\cite{Gab05}]\label{Free products for
  ME}
  If $G_1$ and $G_2$ are measure equivalent to a free group then so is
  $G_1*G_2$.
\end{theorem}

\begin{theorem}[$\mathbf{P_{ME}9}$ in~\cite{Gab05}]\label{Heredity of ME to
  free}
  If $G$ is measure equivalent to a free group and $H\subset G$ is a subgroup
  then $H$ is measure equivalent to a free group.
\end{theorem}
Theorem~\ref{Amalgam of handles} is a special case of a more general result of
Gaboriau.

\begin{theorem}[Corollary 3.19 of~\cite{Gab05}]\label{Gluing a surface along
  one boundary component}
  Consider a countable group $G$ measure equivalent to a free group, and
  $C\subset G$ an infinite cyclic subgroup.  If $\Sigma$ is a compact
  orientable surface of positive genus with a single boundary component then
  $G*_{C=\langle\partial\Sigma\rangle}\pi_1(\Sigma)$ is also measure
  equivalent to a free group.
\end{theorem}
We generalize theorem~\ref{Gluing a surface along one boundary component} to
the case of multiple boundary components.

\begin{corollary}\label{Gluing a surface} Consider a path-connected space $X$
  with $G=\pi_1(X)$ measure equivalent to a free group. Let $\Sigma$ be a
  compact, orientable surface of positive genus with non-empty boundary. Let
  $X'$ be the quotient of $X\sqcup\Sigma$ obtained by identifying the boundary
  curves of $\Sigma$ with loops in $X$ that generate infinite cyclic subgroups
  of $\pi_1(X)$. Then $\pi_1(X')$ is measure equivalent to a free group.
\end{corollary}
\begin{proof}
  By cutting $\Sigma$ along a certain simple closed curve $\gamma$, we can
  decompose it as $\Sigma_1\cup_\gamma\Sigma_2$, where $\Sigma_1$ is planar
  and $\Sigma_2$ is of positive genus and has one boundary component.  The
  space $X'$ acquires a similar decomposition as $X_1\cup_\gamma\Sigma_2$,
  where $X_1$ is obtained from $X$ by amalgamating loops on $X$ with all of
  the boundary curves of $\Sigma_1$ except $\gamma$.  Note that $\Sigma_1$
  deformation retracts onto the graph formed by the boundary circles
  $c_1,\ldots,c_n$ other than $\gamma$, together with a disjoint collection of
  arcs $\alpha_j$ ($j=2,\ldots,n$) connecting $c_1$ to $c_j$.  This
  deformation retraction extends to a deformation retraction of $X_1$ onto the
  union of $X$ and the arcs $\alpha_j$.  It follows from theorem~\ref{Free
  products for ME} that $\pi_1(X_1)\cong\pi_1(X)*F_{n-1}$ is measure
  equivalent to a free group.  Thus
  $\pi_1(X')=\pi_1(X_1)*_{\langle\partial\Sigma_2\rangle}\pi_1(\Sigma_2)$ is
  measure equivalent to a free group, by theorem~\ref{Gluing a surface along
  one boundary component}.
\end{proof}

We are now ready to prove that elementarily free groups are measure equivalent
to free groups.

\subsection{Elementarily free groups}
\begin{theorem}\label{EFGs are ME to free}
  Every elementarily free group is measure equivalent to a free group.
\end{theorem}
\begin{proof}
  By theorem~\ref{EFGs are virtually positive-genus}, it suffices to prove the
  result for positive-genus elementarily free groups.

  At height~$0$, $X_0$ is a one-point union of graphs and hyperbolic surfaces.
  Hyperbolic surface groups are lattices in $\mathrm{PSL}_2(\R)$, so are
  measure equivalent to a free group.  Thus, by theorem~\ref{Free products for
  ME}, $\pi_1(X_0)$ is measure equivalent to a free group.

  At height $h$, assume that $X_h$ is obtained as usual by gluing a surface
  $\Sigma$ to $X_{h-1}$.  By induction, $\pi_1(X_{h-1})$ is measure equivalent
  to a free group.  There are two cases to consider.

  If $\Sigma$ is orientable, then the result is given by corollary~\ref{Gluing
  a surface}.

  If $\Sigma$ is non-orientable, then it has an orientable double cover
  $\Sigma'\rightarrow\Sigma$ of positive genus, with twice the number of
  boundary components.  The amalgam of $\Sigma'$ with two disjoint copies of
  $X_{h-1}$ gives a double cover $X'_h\rightarrow X_h$. Identify a point in
  each copy of $X_{h-1}$ to create a space $Y$.  By proposition~\ref{Free
  products for ME}, $\pi_1(X_{h-1}\vee X_{h-1})=\pi_1(X_{h-1})*\pi_1(X_{h-1})$
  is measure equivalent to a free group.

  We have built $Y$ by gluing the orientable surface of positive genus
  $\Sigma'$ to $X_{h-1}\vee X_{h-1}$, and each boundary component of $\Sigma'$
  defines an element of infinite order in one of the free factors of
  $\pi_1(X_{h-1}\vee X_{h-1})$. It follows by corollary~\ref{Gluing a surface}
  that $\pi_1(Y)$ is measure equivalent to a free group.  Since
  $\pi_1(Y)\cong\pi_1(X'_h)*\Z$, the result follows from theorem~\ref{Heredity
  of ME to free}.
\end{proof}

\section{The case of arbitrary limit groups}
\noindent In the light of theorem~\ref{Heredity of ME to free}, to show that
all limit groups are measure equivalent to free groups it would suffice to
prove that $\omega$-residually free towers are measure equivalent to free
groups.  Even the case of $F_C=F*_{C=Z}\Z^n$, where $C$ is a maximal cyclic
subgroup of $F$ and $Z$ is a direct factor of $\Z^n$, seems non-trivial. The
methods of the proof of theorem~\ref{Gluing a surface along one boundary
component} in~\cite{Gab05} suggest a possible approach.

Let $(X,\mu)$ be a probability measure space, and consider an essentially free
measure-preserving action $\alpha$ of the group $G$ on $X$.  The \emph{orbit
relation} of the action is the equivalence relation given by the orbits of
$G$, and is denoted $\mathcal{R}_\alpha$.  There is a notion of free products
for equivalence relations, motivated by the normal form for free products of
groups.

\begin{definition}
  Consider measured equivalence relations $\mathcal{R},\mathcal{A}$ and
  $\mathcal{B}$ on $X$.  Write $\mathcal{R}=\mathcal{A}*\mathcal{B}$ if:
  \begin{enumerate}
    \item $\mathcal{R}$ is generated by $\mathcal{A}$ and $\mathcal{B}$;
    \item for all $p\ge1$, for almost all $2p$-tuples $(x_i)_{i\in\Z/2p\Z}$ such
      that
      \[
      x_{2j-1}\sim_{\mathcal{A}}x_{2j}\sim_{\mathcal{B}}x_{2j+1}
      \]
      for each $j$, one has $x_i=x_{i+1}$ for some $i$.
  \end{enumerate}
\end{definition}

Gaboriau defines a subgroup $H\subset G$ to be a \emph{measure free factor} if
there exists a free probability-measure-preserving action $\alpha$ of $G$ and
a subrelation $\mathcal{S}$ of $\mathcal{R}_{\alpha}$ so that
\[
\mathcal{R}_{\alpha}=\mathcal{R}_{\alpha|H}*\mathcal{S}.
\]
It follows from the normal form theorem for free products that free factors
are measure free factors.  In~\cite{Gab05}, Gaboriau constructs a non-trivial
example: the boundary circle of a positive-genus orientable surface with one
boundary component generates a measure free factor.  Theorem~\ref{Gluing a
surface along one boundary component} is a special case of:

\begin{theorem}[Theorems 3.13 and 3.17 of~\cite{Gab05}]\label{Gluing along
  measure free factors}
  If $G$ and $G'$ are measure equivalent to free groups, $C\subset G$ and
  $C'\subset G'$ are infinite cyclic subgroups, and $C$ is a measure free
  factor in $G$, then $G*_{C=C'}G'$ is measure equivalent to a free group.
\end{theorem}

At present, the only non-trivial $C$ for which we know that $F_C\me F_2$ is
that given by Gaboriau's example, in which $C$ is generated by a boundary
component of an orientable surface.   It is natural to ask if each maximal
cyclic subgroup of $F$ is a measure free factor.  If this were so, then it
would follow from theorem~\ref{Gluing along measure free factors} that every
$F_C$ is measure equivalent to a free group.  It is also natural to generalize
the question to towers, and ask if each maximal abelian cyclic subgroup of a
tower is a measure free factor.   Again, if so, it would follow that every
limit group is measure equivalent to a free group.

\bibliographystyle{plain}
\bibliography{me}

\end{document}